\documentclass[10pt]{article}
\usepackage{amsmath}
\usepackage{amssymb}
\usepackage{amsthm}
\usepackage{mathrsfs}

\begin{document}

\title{A new estimate for Bochner-Riesz operators at the critical index on the weighted Hardy spaces}
\author{Hua Wang \footnote{E-mail address: wanghua@pku.edu.cn.}\\
\footnotesize{Department of Mathematics, Zhejiang University, Hangzhou 310027, China}}
\date{}
\maketitle

\begin{abstract}
Let $w$ be a Muckenhoupt weight and $H^p_w(\mathbb R^n)$ be the weighted Hardy spaces. In this paper, by using the atomic decomposition of $H^p_w(\mathbb R^n)$, we will show that the Bochner-Riesz operators $T^\delta_R$ are bounded from $H^p_w(\mathbb R^n)$ to the weighted weak Hardy spaces $WH^p_w(\mathbb R^n)$ when $0<p<1$ and $\delta=n/p-(n+1)/2$. This result is new even in the unweighted case.
\end{abstract}

\section{Introduction}

The Bochner-Riesz operators of order $\delta>0$ in $\mathbb R^n$ are defined initially for Schwartz functions in terms of Fourier transforms by
\begin{equation*}
\big(\widehat{T^\delta_R f}\big)(\xi)=\Big(1-\frac{|\xi|^2}{R^2}\Big)^\delta_+\hat{f}(\xi), \quad 0<R<\infty,
\end{equation*}
where $\hat{f}$ denotes the Fourier transform of $f$. The associated maximal Bochner-Riesz operator is defined by
\begin{equation*}
T^\delta_*f(x)=\sup_{R>0}\big|T^\delta_Rf(x)\big|.
\end{equation*}

These operators were first introduced by Bochner \cite{bochner} in connection with summation of multiple Fourier series and played an important role in harmonic analysis. The problem concerning the spherical convergence of Fourier integrals have led to the study of their $L^p$ boundedness. As for their $H^p$ boundedness, Sj\"olin \cite{sjolin} and Stein, Taibleson and Weiss \cite{stein} proved the following theorem (see also [11, page 130]).

\newtheorem*{thmi}{Theorem I}
\begin{thmi}
Suppose that $0<p\le1$ and $\delta>n/p-(n+1)/2$. Then there exists a constant $C>0$ independent of $f$ and $R$ such that
\begin{equation*}
\big\|T^\delta_R(f)\big\|_{H^p}\le C\|f\|_{H^p}.
\end{equation*}
\end{thmi}

In \cite{stein}, the authors also considered weak type estimate for the maximal Bochner-Riesz operator $T^\delta_*$ at the critical index $\delta=n/p-(n+1)/2$ and showed the following inequality is sharp.

\newtheorem*{thmii}{Theorem II}
\begin{thmii}
Suppose that $0<p<1$ and $\delta=n/p-(n+1)/2$. Then there exists a constant $C>0$ independent of $f$ such that
\begin{equation*}
\sup_{\lambda>0}\lambda^p\big|\big\{x\in\mathbb R^n:T^\delta_*f(x)>\lambda\big\}\big|\le C\|f\|^p_{H^p}.
\end{equation*}
\end{thmii}

In 1995, Sato \cite{sato1} studied the weighted case and obtained the following weighted weak type estimate for the maximal Bochner-Riesz operator $T^\delta_*$.

\newtheorem*{thmiii}{Theorem III}
\begin{thmiii}
Let $w\in A_1$$($Muckenhoupt weight class$)$, $0<p<1$ and $\delta=n/p-(n+1)/2$. Then there exists a constant $C>0$ independent of $f$ such that
\begin{equation*}
\sup_{\lambda>0}\lambda^p\cdot w\big(\big\{x\in\mathbb R^n:T^\delta_*f(x)>\lambda\big\}\big)\le C\|f\|_{H^p_w}^p.
\end{equation*}
\end{thmiii}

In \cite{sato2}, Sato also showed that for $n\ge2$, there exists a function $f\in H^1_w\cap L^1$, $w\in A_1$, such that
\begin{equation*}
\limsup_{R\to\infty}\big|T^{{(n-1)}/2}_R f(x)\big|=+\infty \quad \mbox{almost everywhere}.
\end{equation*}
Hence, $\delta=n/p-(n+1)/2$ is indeed the critical index for the weighted case.

In 2006, Lee \cite{lee} considered values of $\delta$ greater than the critical index $n/p-(n+1)/2$ and proved the $H^p_w$--$L^p_w$ boundedness of the maximal operator $T^\delta_*$. Furthermore, by using this strong type estimate of $T^\delta_*$,
Lee \cite{lee} also obtained the $H^p_w$ boundedness of the Bochner-Riesz operator.

The purpose of this article is to discuss the boundedness of Bochner-Riesz operators at the critical index $n/p-(n+1)/2$ on the weighted Hardy spaces. Our main result is stated as follows.

\newtheorem{theorem}{Theorem}[section]
\begin{theorem}
Let $0<p<1$, $\delta=n/p-(n+1)/2$ and $w\in A_1$. Suppose that $n(1/p-1)$ is not a positive integer, then there exists a constant $C>0$ independent of $f$ and $R$ such that
\begin{equation*}
\big\|T^\delta_R(f)\big\|_{WH^p_w}\le C\|f\|_{H^p_w},
\end{equation*}
where $WH^p_w$ denotes the weighted weak Hardy space.
\end{theorem}

In particular, if we take $w$ to be a constant function, then we can get

\newtheorem{cor}[theorem]{Corollary}
\begin{cor}
Let $0<p<1$ and $\delta=n/p-(n+1)/2$. Suppose that $n(1/p-1)$ is not a positive integer, then there exists a constant $C>0$ independent of $f$ and $R$ such that
\begin{equation*}
\big\|T^\delta_R(f)\big\|_{WH^p}\le C\|f\|_{H^p},
\end{equation*}
where $WH^p$ denotes the weak Hardy space.
\end{cor}

\section{Notations and definitions}

The definition of $A_p$ class was first used by Muckenhoupt \cite{muckenhoupt}, Hunt, Muckenhoupt and Wheeden \cite{hunt}, and Coifman and Fefferman \cite{coifman} in the study of weighted
$L^p$ boundedness of Hardy-Littlewood maximal functions and singular integrals. Let $w$ be a nonnegative, locally integrable function defined on $\mathbb R^n$; all cubes are assumed to have their sides parallel to the coordinate axes.
We say that $w\in A_p$, $1<p<\infty$, if
\begin{equation*}
\left(\frac1{|Q|}\int_Q w(x)\,dx\right)\left(\frac1{|Q|}\int_Q w(x)^{-1/{(p-1)}}\,dx\right)^{p-1}\le C \quad\mbox{for every cube}\; Q\subseteq \mathbb
R^n,
\end{equation*}
where $C$ is a positive constant which is independent of the choice of $Q$.

For the case $p=1$, $w\in A_1$, if
\begin{equation*}
\frac1{|Q|}\int_Q w(x)\,dx\le C\cdot\underset{x\in Q}{\mbox{ess\,inf}}\,w(x)\quad\mbox{for every cube}\;Q\subseteq\mathbb R^n.
\end{equation*}

A weight function $w\in A_\infty$ if it satisfies the $A_p$ condition for some $1<p<\infty$. It is well known that if $w\in A_p$ with $1<p<\infty$, then $w\in A_r$ for all $r>p$, and $w\in A_q$ for some $1<q<p$. We thus write $q_w\equiv\inf\{q>1:w\in A_q\}$ to denote the critical index of $w$.

Given a cube $Q$ and $\lambda>0$, $\lambda Q$ denotes the cube with the same center as $Q$ whose side length is $\lambda$ times that of $Q$. $Q=Q(x_0,r)$ denotes the cube centered at $x_0$ with side length $r$. For a weight function $w$ and a measurable set $E$, we denote the Lebesgue measure of $E$ by $|E|$ and set the weighted measure $w(E)=\int_E w(x)\,dx$.

We give the following result that will be used in the sequel.

\newtheorem{lemma}[theorem]{Lemma}
\begin{lemma}[\cite{garcia2}]
Let $w\in A_1$. Then, for any cube $Q$, there exists an absolute constant $C>0$ such that
$$w(2Q)\le C\,w(Q).$$
In general, for any $\lambda>1$, we have
$$w(\lambda Q)\le C\cdot\lambda^{n}w(Q),$$
where $C$ does not depend on $Q$ nor on $\lambda$.
\end{lemma}

\begin{lemma}[\cite{garcia2}]
Let $w\in A_1$. Then there exists a constant $C>0$ such that
\begin{equation*}
C\cdot\frac{|E|}{|Q|}\le\frac{w(E)}{w(Q)}
\end{equation*}
for any measurable subset $E$ of a cube $Q$.
\end{lemma}

Given a weight function $w$ on $\mathbb R^n$, for $0<p<\infty$, we denote by $L^p_w(\mathbb R^n)$ the weighted space of all functions satisfying
\begin{equation*}
\|f\|_{L^p_w}=\left(\int_{\mathbb R^n}|f(x)|^pw(x)\,dx\right)^{1/p}<\infty.
\end{equation*}
We also denote by $WL^p_w(\mathbb R^n)$ the weighted weak $L^p$ space which is formed by all functions satisfying
\begin{equation*}
\|f\|_{WL^p_w}=\sup_{\lambda>0}\lambda\cdot w\big(\big\{x\in\mathbb R^n:|f(x)|>\lambda\big\}\big)^{1/p}<\infty.
\end{equation*}

We write $\mathscr S(\mathbb R^n)$ to denote the Schwartz space of all rapidly decreasing infinitely differentiable functions and $\mathscr S'(\mathbb R^n)$ to denote the space of all tempered distributions, i.e., the topological dual of $\mathscr S(\mathbb R^n)$. For any $0<p\le1$, the weighted Hardy spaces $H^p_w(\mathbb R^n)$ can be defined in terms of maximal functions.
Let $\varphi$ be a function in $\mathscr S(\mathbb R^n)$ satisfying $\int_{\mathbb R^n}\varphi(x)\,dx=1$.
Set
\begin{equation*}
\varphi_t(x)=t^{-n}\varphi(x/t),\quad t>0,\;x\in\mathbb R^n.
\end{equation*}
We will define the maximal function $M^+_\varphi f(x)$ by
\begin{equation*}
M^+_\varphi f(x)=\sup_{t>0}\big|(\varphi_t*f)(x)\big|.
\end{equation*}
Then $H^p_w(\mathbb R^n)$ consists of those tempered distributions $f\in\mathscr S'(\mathbb R^n)$ for which
$M^+_\varphi f\in L^p_w(\mathbb R^n)$ with $\|f\|_{H^p_w}=\|M^+_\varphi f\|_{L^p_w}$. The real-variable theory of weighted Hardy spaces has been investigated by many authors. For example, Garcia-Cuerva \cite{garcia1} studied the atomic decomposition and the dual spaces of $H^p_w$ for $0<p\le1$. The molecular characterization of $H^p_w$ for $0<p\le1$ was given by Lee and Lin \cite{lee2}. We refer the readers to \cite{garcia1,lee2,stomberg} and the references therein for further details.

Let us now turn to the weighted weak Hardy spaces, which are good substitutes for the weighted
Hardy spaces in the study of the boundedness of some operators. The weak $H^p$ spaces have first appeared in the work of Fefferman, Rivi\`ere and Sagher \cite{cfefferman}. The atomic decomposition theory of weak $H^1$ space on $\mathbb R^n$ was given by Fefferman and Soria \cite{rfefferman}. Later, Liu \cite{liu} established the weak $H^p$ spaces on homogeneous
groups. In 2000, Quek and Yang \cite{quek} introduced the weighted weak Hardy spaces and established their atomic decompositions. Let $w\in A_\infty$, $0<p\le1$ and $N=[n(q_w/p-1)]$. Define
\begin{equation*}
\mathscr A_{N,w}=\Big\{\varphi\in\mathscr S(\mathbb R^n):\sup_{x\in\mathbb R^n}\sup_{|\alpha|\le N+1}(1+|x|)^{N+n+1}|D^\alpha\varphi(x)|\le1\Big\},
\end{equation*}
where $\alpha=(\alpha_1,\dots,\alpha_n)\in(\mathbb N\cup\{0\})^n$, $|\alpha|=\alpha_1+\dots+\alpha_n$, and
\begin{equation*}
D^\alpha\varphi=\frac{\partial^{|\alpha|}\varphi}{\partial x^{\alpha_1}_1\cdots\partial x^{\alpha_n}_n}.
\end{equation*}
For $f\in\mathscr S'(\mathbb R^n)$, the nontangential grand maximal function of $f$ is defined by
\begin{equation*}
G_w f(x)=\sup_{\varphi\in\mathscr A_{N,w}}\sup_{|y-x|<t}\big|(\varphi_t*f)(y)\big|.
\end{equation*}
Then we can define the weighted weak Hardy space $WH^p_w(\mathbb R^n)$ by $WH^p_w(\mathbb R^n)=\{f\in\mathscr S'(\mathbb R^n):G_w f\in WL^p_w(\mathbb R^n)\}$ and $\|f\|_{WH^p_w}=\|G_w f\|_{WL^p_w}$.
In order to simplify the computations, we also use another equivalent definition of $WH^p_w(\mathbb R^n)$.
\begin{equation*}
WH^p_w(\mathbb R^n)=\big\{f\in\mathscr S'(\mathbb R^n):G^+_w f\in WL^p_w(\mathbb R^n)\big\},
\end{equation*}
where $G^+_w f$ is called the radial grand maximal function, which is defined by
\begin{equation*}
G^+_w f(x)=\sup_{\varphi\in\mathscr A_{N,w}}\sup_{t>0}\big|(\varphi_t*f)(x)\big|.
\end{equation*}
Moreover, we set $\|f\|_{WH^p_w}=\|G^+_w f\|_{WL^p_w}$.

For $w$ equals to a constant function, we shall denote $WL^p_w(\mathbb R^n)$, $H^p_w(\mathbb R^n)$ and $WH^p_w(\mathbb R^n)$ simply by $WL^p(\mathbb R^n)$, $H^p(\mathbb R^n)$ and $WH^p(\mathbb R^n)$.

Throughout this article $C$ denotes a positive constant, which is independent of the main parameters and not necessarily the same at each occurrence.

\section{The atomic decomposition for weighted Hardy spaces}

In this article, we will use Garcia-Cuerva's atomic decomposition theory for weighted Hardy spaces in \cite{garcia1,stomberg}. We characterize weighted Hardy spaces in terms of atoms in the following way.

Let $0<p\le1\le q\le\infty$ and $p\ne q$ such that $w\in A_q$ with critical index $q_w$. Set [\,$\cdot$\,] the greatest integer function. For $s\in \mathbb Z_+$ satisfying $s\ge N=[n({q_w}/p-1)],$ a real-valued function $a(x)$ is called a ($p,q,s$)-atom centered at $x_0$ with respect to $w$ (or a $w$-($p,q,s$)-atom centered at $x_0$) if the following conditions are satisfied:

(a) $a\in L^q_w(\mathbb R^n)$ and is supported in a cube $Q$ centered at $x_0$;

(b) $\|a\|_{L^q_w}\le w(Q)^{1/q-1/p}$;

(c) $\int_{\mathbb R^n}a(x)x^\alpha\,dx=0$ for every multi-index $\alpha$ with $|\alpha|\le s$.

\begin{theorem}
Let $0<p\le1\le q\le\infty$ and $p\ne q$ such that $w\in A_q$ with critical index $q_w$. For each $f\in H^p_w(\mathbb R^n)$, there exist a sequence \{$a_j$\} of $w$-$(p,q,N)$-atoms and a sequence \{$\lambda_j$\} of real numbers with $\sum_j|\lambda_j|^p\le C\|f\|^p_{H^p_w}$ such that $f=\sum_j\lambda_j a_j$ both in the sense of distributions and in the $H^p_w$ norm.
\end{theorem}

\section{Proof of Theorem 1.1}

The Bochner-Riesz operators can be expressed as convolution operators
\begin{equation*}
T^\delta_Rf(x)=(\phi_{1/R}*f)(x),
\end{equation*}
where $\phi(x)=[(1-|\cdot|^2)^\delta_+]\mbox{\textasciicircum}(x)$ and $\phi_{1/R}(x)=R^n\phi(Rx)$. It is well known that the kernel $\phi$ can be represented as (see \cite{lu2,stein2})
\begin{equation*}
\phi(x)=\pi^{-\delta}\Gamma(\delta+1)|x|^{-(\frac n2+\delta)}J_{\frac n2+\delta}(2\pi|x|),
\end{equation*}
where $J_\mu(t)$ is the Bessel function
\begin{equation*}
J_\mu(t)=\frac{(\frac{t}{2})^\mu}{\Gamma(\mu+\frac12)\Gamma(\frac12)}\int_{-1}^1e^{its}(1-s^2)^{\mu-\frac12}\,ds.
\end{equation*}
The following kernel estimates of these convolution operators are well known. For its proof, we refer the readers to \cite{sato1}. See also [11, page 121].

\begin{lemma}
Let $0<p<1$ and $\delta=n/p-(n+1)/2$. Then the kernel $\phi$ satisfies the inequality
\begin{equation*}
\sup_{x\in\mathbb R^n}(1+|x|)^{n/p}\big|D^\alpha\phi(x)\big|\le C\quad\mbox{for all multi-indices}\;\alpha.
\end{equation*}
\end{lemma}

In order to prove our main theorem, we shall need the following superposition principle on weighted weak type estimates.

\begin{lemma}
Let $w\in A_1$ and $0<p<1$.
If a sequence of measurable functions $\{f_j\}$ satisfy
\begin{equation*}
w\big(\big\{x\in\mathbb R^n:|f_j(x)|>\alpha\big\}\big)\le \alpha^{-p} \quad\mbox{for all}\;\, j\in\mathbb Z
\end{equation*}
and
\begin{equation*}
\sum_{j\in\mathbb Z}|\lambda_j|^p\le1,
\end{equation*}
then we have that $\sum_j\lambda_jf_j(x)$ is absolutely convergent almost everywhere and
\begin{equation*}
w\Big(\Big\{x\in\mathbb R^n:\big|\sum_j\lambda_jf_j(x)\big|>\alpha\Big\}\Big)\le\frac{2-p}{1-p}\cdot\alpha^{-p}.
\end{equation*}
\end{lemma}
\begin{proof}
The proof of this lemma is similar to the corresponding result for the unweighted case which can be found in \cite{stein}. See also [11, page 123].
\end{proof}

We now prove the following estimate which plays a key role in the proof of our main result.

\begin{lemma}
Let $0<p<1<q<\infty$, $\delta=n/p-(n+1)/2$ and $w\in A_1$. Then for any given $w$-$(p,q,N)$-atom $a$ with support contained in $Q=Q(x_0,r)$, we have
\begin{equation*}
T^\delta_*(a)(y)\le C\cdot\frac{r^{n/p}}{|y-x_0|^{n/p}}w(Q)^{-1/p}\quad \mbox{whenever}\; \;|y-x_0|>\sqrt{n}\,r.
\end{equation*}
Here, and in what follows, we always denote $N=[n({q_w}/p-1)]=[n(1/p-1)]$.
\end{lemma}
\begin{proof}
Actually, this lemma was essentially contained in \cite{lee}. For the sake of completeness, we give its proof here. For any $\varepsilon>0$, we write
\begin{equation*}
(\phi_\varepsilon*a)(y)=\varepsilon^{-n}\int_{Q(x_0,r)}\phi\Big(\frac{y-z}{\varepsilon}\Big)a(z)\,dz.
\end{equation*}
Let us now consider the following two cases.

$(i)$ $0<\varepsilon\le r$. Note that $\delta=n/p-(n+1)/2$, then by Lemma 4.1, we have
\begin{equation*}
\big|(\phi_\varepsilon*a)(y)\big|\le C\cdot\varepsilon^{n/p-n}\int_{Q(x_0,r)}\frac{|a(z)|}{|y-z|^{n/p}}\,dz.
\end{equation*}
By our assumption, when $z\in Q(x_0,r)$, then we can easily get $|y-z|\ge|y-x_0|-|z-x_0|\ge\frac{|y-x_0|}{2}$. Observe that $0<p<1$, then $n/p-n>0$. Hence
\begin{equation}
\big|(\phi_\varepsilon*a)(y)\big|\le C\cdot\frac{r^{n/{p}-n}}{|y-x_0|^{n/{p}}}\int_{Q(x_0,r)}|a(z)|\,dz.
\end{equation}
Denote the conjugate exponent of $q>1$ by $q'={q}/(q-1).$ Using H\"older's inequality, the $A_q$ condition and the size condition of atom $a$, we can get
\begin{align}
\int_{Q(x_0,r)}|a(z)|\,dz&\le\bigg(\int_{Q(x_0,r)}\big|a(z)\big|^{q}w(z)\,dz\bigg)^{1/q}
\bigg(\int_{Q(x_0,r)}\big(w(z)^{-1/q}\big)^{q'}\,dz\bigg)^{1/{q'}}\notag\\
&\le C\cdot\|a\|_{L^{q}_{w}}|Q(x_0,r)|w(Q)^{-1/q}\notag\\
&\le C\cdot|Q(x_0,r)|w(Q)^{-1/p}.
\end{align}
Substituting the above inequality (2) into (1), we thus obtain
\begin{equation}
\big|(\phi_\varepsilon*a)(y)\big|\le C\cdot\frac{r^{n/{p}}}{|y-x_0|^{n/{p}}}w(Q)^{-1/p}.
\end{equation}

$(ii)$ $\varepsilon>r$. It is easy to see that the choice of $N$ in this lemma implies $\frac{n}{n+N+1}<p\le\frac{n}{n+N}$. Using the vanishing moment condition of atom $a$, Taylor's theorem and Lemma 4.1, we deduce
\begin{equation*}
\begin{split}
\big|(\phi_\varepsilon*a)(y)\big|&=\varepsilon^{-n}\bigg|\int_{Q(x_0,r)}
\Big[\phi\Big(\frac{y-z}{\varepsilon}\Big)-\sum_{|\gamma|\le N}\frac{D^\gamma\phi(\frac{y-x_0}{\varepsilon})}{\gamma!}\Big(\frac{z-x_0}{\varepsilon}\Big)^\gamma\Big] a(z)\,dz\bigg|\\
&\le\varepsilon^{-n}\cdot\Big(\frac{{\sqrt n}r}{2\varepsilon}\Big)^{N+1}
\int_{Q(x_0,r)}\sum_{|\gamma|=N+1}\Big|\frac{D^\gamma\phi(\frac{y-x_0-\theta (z-x_0)}{\varepsilon})}{\gamma!}\Big|\big|a(z)\big|\,dz\\
&\le C\cdot\frac{r^{N+1}}{\varepsilon^{n+N+1}}\int_{Q(x_0,r)}\Big|\frac{y-x_0-\theta (z-x_0)}{\varepsilon}\Big|^{-n/{p}}\big|a(z)\big|\,dz,
\end{split}
\end{equation*}
where $0<\theta<1$. As in the first case $(i)$, we have $|y-x_0|\ge2|z-x_0|$, which implies $\big|y-x_0-\theta (z-x_0)\big|\ge|y-x_0|-|z-x_0|\ge\frac{|y-x_0|}{2}$. This together with the inequality (2) yields
\begin{equation*}
\big|(\phi_\varepsilon*a)(y)\big|\le C\cdot\frac{r^{n+N+1}}{\varepsilon^{n+N+1-n/{p}}}\frac{1}{|y-x_0|^{n/{p}}}w(Q)^{-1/p}.
\end{equation*}
Observe that $n+N+1-n/{p}>0$, then for $\varepsilon>r$, we have $\varepsilon^{n+N+1-n/{p}}>r^{n+N+1-n/{p}}$. Consequently
\begin{equation}
\big|(\phi_\varepsilon*a)(y)\big|\le C\cdot\frac{r^{n/{p}}}{|y-x_0|^{n/{p}}}w(Q)^{-1/p}.
\end{equation}
Summarizing the estimates (3) and (4) derived above and taking the supremum over all $\varepsilon>0$, we obtain the desired inequality.
\end{proof}

We are now in a position to give the proof of Theorem 1.1.

\begin{proof}
By the atomic decomposition theorem (Theorem 3.1) and Lemma 4.2, it suffices to show that for any $w$-$(p,q,N)$-atom $a$, there exists a constant $C>0$ independent of $a$ such that $\|G^+_w(T^\delta_Ra)\|_{WL^p_w}\le C$.

Let $a$ be a $w$-$(p,q,N)$-atom with $supp\, a\subseteq Q=Q(x_0,r)$, $0<p<1<q<\infty$, and let $Q^*=(4\sqrt n)Q$. For any given $\lambda>0$, we write
\begin{equation*}
\begin{split}
&\lambda^p\cdot w\big(\big\{x\in\mathbb R^n:|G^+_w(T^\delta_Ra)(x)|>\lambda\big\}\big)\\
\le\,&\lambda^p\cdot w\big(\big\{x\in Q^*:|G^+_w(T^\delta_Ra)(x)|>\lambda\big\}\big)+\lambda^p\cdot w\big(\big\{x\in(Q^*)^c:|G^+_w(T^\delta_Ra)(x)|>\lambda\big\}\big)\\
=\,&I_1+I_2.
\end{split}
\end{equation*}
Let us deal with the term $I_1$ first. Applying Chebyshev's inequality and H\"older's inequality, we obtain
\begin{equation*}
\begin{split}
I_1&\le\int_{Q^*}\big|G^+_w(T^\delta_Ra)(x)\big|^pw(x)\,dx\\
&\le\bigg(\int_{Q^*}\big|G^+_w(T^\delta_Ra)(x)\big|^qw(x)\,dx\bigg)^{p/q}\bigg(\int_{Q^*}w(x)\,dx\bigg)^{1-p/q}\\
&\le \big\|G^+_w(T^\delta_Ra)\big\|^p_{L^q_w}w(Q^*)^{1-p/q}.
\end{split}
\end{equation*}
Since $0<p<1$ and $\delta=n/p-(n+1)/2$, then $\delta>(n-1)/2$. In this case, it is well known that (see \cite{lu2,stein2})
\begin{equation*}
T^\delta_* a(x)\le C\cdot M a(x),
\end{equation*}
where $M$ denotes the Hardy-Littlewood maximal operator. In addition, we can easily see that for any function $f$, the following inequality holds
\begin{equation*}
G^+_w f(x)\le C\cdot M f(x).
\end{equation*}
Since $w\in A_1$, then $w\in A_q$ for any $1<q<\infty$. Hence, it follows from the $L^q_w$ boundedness of $M$ and Lemma 2.1 that
\begin{align}
I_1&\le C\cdot\big\|T^\delta_*(a)\big\|^p_{L^q_w}w(Q^*)^{1-p/q}\notag\\
&\le C\cdot\|a\|^p_{L^q_w}w(Q)^{1-p/q}\notag\\
&\le C.
\end{align}

We now turn to estimate the other term $I_2$. First, we claim that for every multi-index $\gamma$ with $|\gamma|\le N$ and for any $0<R<\infty$, the following equality holds
\begin{equation}
\int_{\mathbb R^n}T^\delta_Ra(y)y^{\gamma}\,dy=0.
\end{equation}
In fact, by using the Leibnitz formula and vanishing moment condition of atom $a$, we thus have
\begin{equation*}
\begin{split}
\int_{\mathbb R^n}T^\delta_Ra(y)y^{\gamma}\,dy&=(T^\delta_Ra(y)y^{\gamma}){\mbox{\textasciicircum}}(0)\\
&=C\cdot D^{\gamma}(\widehat{T^\delta_Ra})(0)\\
&=C\cdot D^{\gamma}(\widehat{\phi_{1/R}}\cdot\widehat a)(0)\\
&=C\cdot\sum_{|\alpha|+|\beta|=|\gamma|}(D^\alpha\widehat{\phi_{1/R}})(0)\cdot(D^\beta\widehat a)(0)\\
&=0,
\end{split}
\end{equation*}
which is the desired conclusion. Then by the equality (6) derived above, we can split the following expression into three parts and estimate each term respectively.
\begin{equation*}
\begin{split}
\big|\varphi_t*(T^\delta_Ra)(x)\big|=&\bigg|\int_{\mathbb R^n}t^{-n}\Big[\varphi\Big(\frac{x-y}{t}\Big)-\sum_{|\gamma|\le N}\frac{D^\gamma\varphi(\frac{x-x_0}{t})}{\gamma!}\Big(\frac{y-x_0}{t}\Big)^\gamma\Big]T^\delta_Ra(y)\,dy\bigg|\\
\le&\,t^{-n}\bigg|\int_{|y-x_0|\le{\sqrt n} r}\Big[\varphi\Big(\frac{x-y}{t}\Big)-\sum_{|\gamma|\le N}\frac{D^\gamma\varphi(\frac{x-x_0}{t})}{\gamma!}\Big(\frac{y-x_0}{t}\Big)^\gamma\Big]T^\delta_Ra(y)\,dy\bigg|\\
&+t^{-n}\bigg|\int_{{\sqrt n} r<|y-x_0|\le\frac{|x-x_0|}{2}}\Big[\varphi\Big(\frac{x-y}{t}\Big)-\sum_{|\gamma|\le N}\frac{D^\gamma\varphi(\frac{x-x_0}{t})}{\gamma!}\Big(\frac{y-x_0}{t}\Big)^\gamma\Big]T^\delta_Ra(y)\,dy\bigg|\\
&+t^{-n}\bigg|\int_{|y-x_0|>\frac{|x-x_0|}{2}}\Big[\varphi\Big(\frac{x-y}{t}\Big)-\sum_{|\gamma|\le N}\frac{D^\gamma\varphi(\frac{x-x_0}{t})}{\gamma!}\Big(\frac{y-x_0}{t}\Big)^\gamma\Big]T^\delta_Ra(y)\,dy\bigg|\\
=&\,J_1+J_2+J_3.
\end{split}
\end{equation*}
For the term $J_1$, by using Taylor's theorem, we then obtain
\begin{equation*}
J_1\le t^{-n}\cdot\Big(\frac{{\sqrt n}r}{t}\Big)^{N+1}
\int_{|y-x_0|\le\sqrt{n}r}\sum_{|\gamma|=N+1}\Big|\frac{D^\gamma\varphi(\frac{x-x_0-\theta' (y-x_0)}{t})}{\gamma!}\Big|\big|T^\delta_*a(y)\big|\,dy,
\end{equation*}
where $0<\theta'<1$. For any $x\in(Q^*)^c$ and $|y-x_0|\le\sqrt{n}r$, as before, we have $\big|x-x_0-\theta'(y-x_0)\big|\ge\frac{|x-x_0|}{2}$. Observe that $\varphi\in\mathscr A_{N,w}$, then it follows from H\"older's inequality and the size condition of atom $a$ that
\begin{equation*}
\begin{split}
J_1&\le C\cdot\frac{r^{N+1}}{t^{n+N+1}}\int_{|y-x_0|\le\sqrt{n}r}\Big|\frac{x-x_0-\theta' (y-x_0)}{t}\Big|^{-n-N-1}\big|T^\delta_*a(y)\big|\,dy\\
&\le C\cdot\frac{r^{N+1}}{|x-x_0|^{n+N+1}}\big\|T^\delta_*(a)\big\|_{L^q_w}\big|Q(x_0,2\sqrt{n} r)\big|
w\big(Q(x_0,2\sqrt{n} r)\big)^{-1/q}\\
&\le C\cdot\frac{r^{N+1}}{|x-x_0|^{n+N+1}}\|a\|_{L^q_w}\big|Q(x_0,r)\big|
w(Q)^{-1/q}\\
&\le C\cdot\frac{r^{n+N+1}}{|x-x_0|^{n+N+1}}w(Q)^{-1/p}.
\end{split}
\end{equation*}
On the other hand, for the term $J_2$, notice that $|y-x_0|>\sqrt{n}r$ and $\varphi\in\mathscr A_{N,w}$, then by Taylor's theorem and Lemma 4.3, we deduce
\begin{equation*}
\begin{split}
J_2&\le t^{-n}\int_{\sqrt{n}r<|y-x_0|\le\frac{|x-x_0|}{2}}
\sum_{|\gamma|=N+1}\Big|\frac{D^\gamma\varphi(\frac{x-x_0-\theta' (y-x_0)}{t})}{\gamma!}\Big|\Big|\frac{y-x_0}{t}\Big|^{N+1}\big|T^\delta_*a(y)\big|\,dy\\
&\le C\cdot t^{-n}\int_{\sqrt{n}r<|y-x_0|\le\frac{|x-x_0|}{2}}\frac{t^{n+N+1}}{|x-x_0|^{n+N+1}}\cdot\frac{|y-x_0|^{N+1}}{t^{N+1}}
\cdot\frac{r^{n/p}}{|y-x_0|^{n/p}}w(Q)^{-1/p}\,dy\\
&\le C\cdot\frac{r^{n/p}}{|x-x_0|^{n+N+1}}w(Q)^{-1/p}
\int_{\sqrt{n}r<|y-x_0|\le\frac{|x-x_0|}{2}}\big|y-x_0\big|^{N+1-n/p}\,dy.
\end{split}
\end{equation*}
If we rewrite this integral in polar coordinates and note the fact that $-1<n+N-n/p\le0$, then we can get
\begin{equation*}
\begin{split}
\int_{\sqrt{n}r<|y-x_0|\le\frac{|x-x_0|}{2}}\big|y-x_0\big|^{N+1-n/p}\,dy&\le
C\cdot\int_{\sqrt{n}r}^{\frac{|x-x_0|}{2}}\rho^{n+N-n/p}\,d\rho\\
&\le C\cdot\big|x-x_0\big|^{n+N+1-n/p}.
\end{split}
\end{equation*}
So we have
\begin{equation*}
J_2\le C\cdot\frac{r^{n/{p}}}{|x-x_0|^{n/{p}}}w(Q)^{-1/p}.
\end{equation*}
For the last term $J_3$, when $x\in (Q^*)^c$, then it is easy to see that $|y-x_0|>\frac{|x-x_0|}{2}>\sqrt{n}r$. Therefore, it follows from Lemma 4.3 and the fact $\varphi\in\mathscr A_{N,w}$ that
\begin{equation*}
\begin{split}
J_3&\le t^{-n}\int_{|y-x_0|>\frac{|x-x_0|}{2}}\bigg\{\Big|\varphi\Big(\frac{x-y}{t}\Big)\Big|+\sum_{j=0}^N\sum_{|\gamma|=j}
\Big|\frac{D^\gamma\varphi(\frac{x-x_0}{t})}{\gamma!}\Big|\Big|\frac{y-x_0}{t}\Big|^j\bigg\}
\big|T^\delta_*a(y)\big|\,dy\\
&\le C\int_{|y-x_0|>\frac{|x-x_0|}{2}}\bigg\{\big|\varphi_t(x-y)\big|+\sum_{j=0}^N\frac{|y-x_0|^j}{|x-x_0|^{n+j}}
\bigg\}\cdot\frac{r^{n/{p}}}{|y-x_0|^{n/{p}}}w(Q)^{-1/p}\,dy\\
&\le C\cdot w(Q)^{-1/p}\bigg(\frac{r^{n/{p}}}{|x-x_0|^{n/{p}}}\|\varphi_t\|_{L^1}+
\sum_{j=0}^N\frac{r^{n/{p}}}{|x-x_0|^{n+j}}\int_{|y-x_0|>\frac{|x-x_0|}{2}}\frac{dy}{|y-x_0|^{n/p-j}}\bigg).
\end{split}
\end{equation*}
By our assumption (say $n(1/p-1)$ is not a positive integer), we know that $n/p-n>N=[n(1/p-1)]$. Thus we have $n/p-n-j>0$ for any $0\le j\le N$. Making use of the polar coordinates again, we obtain
\begin{equation*}
\begin{split}
\int_{|y-x_0|>\frac{|x-x_0|}{2}}\frac{dy}{|y-x_0|^{n/p-j}}&\le C\cdot
\int_{\frac{|x-x_0|}{2}}^\infty\frac{1}{\rho^{n/p-n-j+1}}\,d\rho\\
&\le C\cdot\frac{1}{|x-x_0|^{n/p-n-j}}.
\end{split}
\end{equation*}
Hence
\begin{equation*}
J_3\le C\cdot\frac{r^{n/{p}}}{|x-x_0|^{n/{p}}}w(Q)^{-1/p}.
\end{equation*}

Therefore, combining the above estimates for $J_1$, $J_2$ and $J_3$, and then taking the supremum over all $t>0$ and all $\varphi\in\mathscr A_{N,w}$, we finally obtain that for any $x\in(Q^*)^c$,
\begin{align}
\big|G^+_w(T^\delta_Ra)(x)\big|&\le C\cdot\left(\frac{r^{n+N+1}}{|x-x_0|^{n+N+1}}w(Q)^{-1/p}+\frac{r^{n/{p}}}{|x-x_0|^{n/{p}}}w(Q)^{-1/p}\right)\notag\\
&\le C\cdot\frac{r^{n/{p}}}{|x-x_0|^{n/{p}}}w(Q)^{-1/p}.
\end{align}

Without loss of generality, we may assume that the absolute constant $C$ appearing in the above inequality is bigger than one. For any fixed $\lambda>0$, we are going to consider two cases. If $\lambda\ge C\cdot\frac{r^{n/{p}}}{{(2\sqrt nr)}^{n/{p}}}w(Q)^{-1/p}$,
then for any $x\in(Q^*)^c$, we have $|x-x_0|\ge2\sqrt nr$. Thus in view of (7),
we can easily see that
$\big\{x\in(Q^*)^c:|G^+_w(T^\delta_Ra)(x)|>\lambda\big\}=\O$. Hence, in this case, the inequality
$$I_2\le C$$
holds trivially. Now suppose that $\lambda<C\cdot\frac{r^{n/{p}}}{{(2\sqrt nr)}^{n/{p}}}w(Q)^{-1/p}$. If we set 
\begin{equation*}
R=\left(\frac{C^p\cdot|Q|}{\lambda^p\cdot w(Q)}\right)^{1/n},
\end{equation*}
then it is not difficult to verify that $R\ge2\sqrt nr\ge r$ and
\begin{equation*}
\big\{x\in(Q^*)^c:|G^+_w(T^\delta_Ra)(x)|>\lambda\big\}\subseteq\big\{x\in\mathbb R^n:|x-x_0|<R\big\}\subseteq Q(x_0,2R).
\end{equation*}
Since $w\in A_1$, then by Lemma 2.2, we can get (below, $\widetilde C$ is an absolute constant)
\begin{equation*}
\widetilde C \cdot\frac{|Q(x_0,r)|}{|Q(x_0,R)|}\le\frac{w(Q(x_0,r))}{w(Q(x_0,R))},
\end{equation*}
which is equivalent to
\begin{equation*}
\begin{split}
w\big(Q(x_0,R)\big)&\le\frac{R^n\cdot w(Q)}{\widetilde C\cdot|Q|}\\
&\le\frac{C^p}{\widetilde C\cdot\lambda^p}.
\end{split}
\end{equation*}
Therefore 
\begin{align}
I_2&\le\lambda^p\cdot w\big(Q(x_0,2R)\big)\notag\\
&\le C^p\cdot{\widetilde C}^{-1}.
\end{align}
Combining the above estimate (8) with (5) and then taking the supremum over all $\lambda>0$, we complete the proof of Theorem 1.1.
\end{proof}

\newtheorem*{rek}{Remark}
\begin{rek}
In the end, we would like to point out that the conclusions of \mbox{\upshape Theorem 1.1} and \mbox{\upshape Corollary 1.2} are not generally true for the maximal Bochner-Riesz operators $T^\delta_*$, because we can not deduce that $T^\delta_*$ satisfy the previous vanishing moment condition $(6)$.
\end{rek}

\end{document}